\documentclass[a4paper, 10pt]{article}

\usepackage[latin1]{inputenc}
\usepackage[T1]{fontenc}
\usepackage{amsmath,amsthm,amssymb}
\usepackage[english]{babel}
\usepackage[all]{xy}
\usepackage{fancyhdr}
\usepackage{geometry}
\title{On maximal diagonalizable Lie subalgebras of the first Hochschild cohomology}
\author{Patrick Le Meur}
\date{\today}

\theoremstyle{definition}
\newtheorem{definition}{Definition}[section]
\theoremstyle{plain}
\newtheorem{Thm}[]{Theorem}
\newtheorem{thm}[]{Theorem}

\newtheorem{lem}[definition]{Lemma}
\newtheorem{prop}[definition]{Proposition}
\theoremstyle{remark}
\newtheorem{ex}[definition]{Example}

\def\c#1{\mathcal{#1}}
\def\hom{\text{\normalfont\textsf{Hom}}}
\def\im{\text{\normalfont\textsf{Im}}}
\def\ker{\text{\normalfont\textsf{Ker}}}
\def\dim{\text{\normalfont\textsf{dim}}}
\def\carac{\text{\normalfont\textsf{char}}}
\def\ext{\text{\normalfont\textsf{Ext}}}
\def\hh{\text{\normalfont\textsf{HH}}}
\def\e{\varepsilon}
\def\square{\blacksquare}

\def\fd{finite dimensional }
\def\kalg{$k$-algebra }
\def\bfda{basic finite dimensional $k$-algebra }
\def\bcfda{basic connected finite dimensional $k$-algebra }

\makeatletter
\renewcommand\section{\@startsection {section}{1}{0mm}%
{-3.5ex \@plus -1ex \@minus -.2ex}
                                   {0.5ex \@plus.2ex}%
                                   {\normalfont\large\bfseries}}

\renewcommand\subsection{\@startsection {subsection}{2}{5mm}%
  {-3.5ex \@plus -1ex \@minus -.2ex}
  {0.5ex \@plus.2ex}
  {\normalfont\normalsize\itshape}
}

\makeatother

\begin{document}
\maketitle
\abstract{Let $A$ be a basic connected finite dimensional algebra
  over an algebraically closed field, with ordinary quiver without
  oriented cycles. Given a presentation of $A$ by quiver
  and admissible relations, Assem and de la Pe\~na have constructed an
  embedding of the space of additive characters of the fundamental
  group of the presentation into the first Hochschild cohomology group
  of $A$. We
  compare the embeddings given by the different presentations of
  $A$. In some situations, we characterise the images of these
  embeddings in terms of (maximal) diagonalizable subalgebras of the
  first Hochschild cohomology group (endowed with its Lie algebra structure).}
\section*{Introduction}
Let $A$ be a finite dimensional algebra over an algebraically closed
field $k$. The representation theory of $A$ deals with the study of (right)
$A$-modules. So we assume that $A$ is basic and connected and it admits
presentations $A\simeq kQ/I$ by its (unique) ordinary quiver $Q$ and
an ideal $I$ of admissible relations. In the eighties, Martinez-Villa and de la Pe\~na
introduced the fundamental group $\pi_1(Q,I)$ of $(Q,I)$ (\cite{martinezvilla_delapena}). Like in topology, this group is
defined using an equivalence relation $\sim_I$ (called the homotopy
relation) on the set of unoriented paths in $Q$. This group is part of the so-called
covering techniques initiated in \cite{bongartz_gabriel,riedtmann}. In particular, it has led to the definition
of simple connectedness and strong simple connectedness for an algebra
(\cite{assem_skowronski,skowronski2}). Also, it has proved to be a very useful tool
in representation theory. For example, it is proved in
\cite{skowronski3} that any domestic self-injective algebra admitting a Galois
covering by a strongly simply connected locally bounded $k$-category
is of quasi-tilted type. Note that in general, different presentations
$A\simeq kQ/I$ and $A\simeq kQ/J$ may lead to non-isomorphic groups
$\pi_1(Q,I)$ and $\pi_1(Q,J)$.

The fundamental group $\pi_1(Q,I)$ behaves much like the fundamental
group of a topological space. For example, given a presentation
$\nu\colon kQ\twoheadrightarrow A$ (with kernel $I$), Assem and de la Pe\~na
have defined an injective group homomorphism
$\theta_{\nu}\colon\hom(\pi_1(Q,I),k^+)\hookrightarrow \hh^1(A)$. Here
$\hh^1(A)$ is the
first Hochschild cohomology group $\ext_{A^{op}\otimes A}^1(A,A)$ (\cite{hochschild}). This result is to be compared with  the classical isomorphism $\hom(\pi_1(X),\mathbb{Z})\xrightarrow{\sim} H^1(X;\mathbb{Z})$
relating the singular cohomology to the fundamental group of a path connected topological space $X$.
It is known from \cite{gerstenhaber} that
$\hh^1(A)$ has a structure of Lie algebra, isomorphic to
the Lie algebra of derivations of $A$ (with the commutator as Lie bracket)
factored out by the ideal of inner
derivations. With this presentation of $\hh^1(A)$, the derivations
that lie in the image of $\theta_{\nu}$, have been characterized in
terms of diagonalizable derivations (\cite{farkas_green_marcos}, see also \cite{delapena_saorin}). 

The aim of this text is
to characterise
 maximal diagonalisable Lie subalgebras of $\hh^1(A)$ using the
 subspaces $\im(\theta_{\nu})$ associated to the
different  presentations $\nu$ of $A$. Recall that, given a Lie algebra,
the maximal diagonalizable (for the adjoint representation)
subalgebras are related to Cartan subalgebras.

On the one hand, one can define a diagonalizability for elements in
$\hh^1(A)$ using the above notion of diagonalizable derivations. Also
one can define the diagonalizability of a subset of $\hh^1(A)$ (as the
simultaneous diagonalizability of its elements). It appears that
$\im(\theta_{\nu})$ is diagonalizable, and that any diagonalizable
subset of $\hh^1(A)$ is contained in $\im(\theta_{\nu})$ for some
presentation $\nu\colon kQ\twoheadrightarrow A$.

On the other hand, given two presentations $\nu\colon
kQ\twoheadrightarrow A$ and $\mu\colon kQ\twoheadrightarrow A$ with
kernel $I$ and $J$ respectively, it is not easy to compare the groups
$\pi_1(Q,I)$ and $\pi_1(Q,J)$ (and therefore
$\theta_{\nu}$ and $\theta_{\mu}$). In some cases, this is possible,
however. For example, assume that $(\alpha,u)$ is a bypass in $Q$
(that is, $\alpha$ is an arrow and $u$ is an oriented path which is
parallel to $\alpha$ and distinct from $\alpha$), that $\tau\in k$,
and that $J=\varphi_{\alpha,u,\tau}(I)$. Here
$\varphi_{\alpha,u,\tau}\colon kQ\xrightarrow{\sim}kQ$ is the
automorphism, called a \emph{transvection}, which maps $\alpha$ to
$\alpha+\tau u$, and which fixes any other arrow (\cite{lemeur}). 
In such a situation, if $\alpha\sim_Iu$ (or $\alpha\sim_Ju$), then
there is a natural surjective group homomorphism 
$\pi_1(Q,J)\twoheadrightarrow\pi_1(Q,I)$ (or
$\pi_1(Q,I)\twoheadrightarrow\pi_1(Q,J)$, respectively); if
$\alpha\sim_Iu$ and $\alpha\sim_Ju$ then $\pi_1(Q,I)$ equals
$\pi_1(Q,J)$ and the natural homomorphisms are the identity maps; and
if  $\alpha\not\sim_Iu$ and $\alpha\not\sim_Ju$, then $I=J$ and
$\pi_1(Q,I)=\pi_1(Q,J)$. 
In each of these cases, we shall see that there is a simple relation between
$\theta_{\nu}$ and $\theta_{\mu}$.

In order to formulate our main result, we use the quiver $\Gamma$ of the homotopy
relations of the presentations of $A$ (\cite{lemeur}). Its set of vertices is the set
of the homotopy relations $\sim_{\ker(\nu)}$ associated to all the
presentations $\nu\colon kQ\twoheadrightarrow A$. Also, there is an
arrow $\sim_I\to \sim_J$ if there exists a transvection
$\varphi_{\alpha,u,\tau}$ such that $J=\varphi_{\alpha,u,\tau}(I)$ and
such that the natural surjective group homomorphism is a non isomorphism
$\pi_1(Q,I)\twoheadrightarrow\pi_1(Q,J)$. The quiver $\Gamma$ has been
introduced in order to find conditions under which an algebra admits a
universal Galois covering. This existence is related to the existence
of a unique source (that is, a vertex which is the target of no
arrow). Actually, under one of the two following conditions, $\Gamma$
does have a unique source
(\cite[Prop. 2.11]{lemeur2} and \cite[Cor. 4.4]{lemeur3}):
\begin{enumerate}
\item[($H_1$)] $Q$ has no double bypass and $k$ has characteristic
  zero (a \emph{double bypass} is a $4$-tuple $(\alpha,u,\beta,v)$
  where $(\alpha,u)$ and $(\beta,v)$ are bypasses such that the arrow
  $\beta$ appears in the path $u$).
\item[($H_2$)] $A$ is monomial  (that is, $A\simeq kQ/I_0$ with $I_0$ an ideal
  generated by a set of paths) and $Q$ has no multiple arrows.
\end{enumerate}

Using these results, we prove the main theorem of the text.
\begin{Thm}
\label{thm1}
Assume that at least one the two hypotheses ($H_1$) or ($H_2$) is
satisfied. Then:
\begin{enumerate}
\item[($i$)] The maximal diagonalizable subalgebras of $\hh^1(A)$ are
  exactly the subalgebras of the form $\im(\theta_{\nu})$ where $\nu\colon
  kQ\twoheadrightarrow A$ is a presentation such that
  $\sim_{\ker(\nu)}$ is the unique source of $\Gamma$.
\item[($ii$)] If $\mathcal{G},\mathcal{G}'$ are two such subalgebras
  of $\hh^1(A)$, then there exists an algebra automorphism $\psi\colon
  A\xrightarrow{\sim} A$ inducing a Lie algebra automorphism
  $\psi_*\colon \hh^1(A)\xrightarrow{\sim}\hh^1(A)$ such that $\mathcal{G}'=\psi_*(\mathcal{G})$.
\end{enumerate}
\end{Thm}

Note that the Lie algebra $\hh^1(A)$ has already been studied (see
\cite{guilasensio_saorin,strametz}, for instance).

The text is organised as follows. In Section~$1$ we recall all
the definitions we will need and prove some useful lemmas. In
Section~$2$, we introduce the notion of diagonalizability in
$\hh^1(A)$. In particular, we will prove that a subset of $\hh^1(A)$ is
diagonalizable is
and only if it is contained in $\im(\theta_{\nu})$ for some 
 presentation $\nu\colon kQ\twoheadrightarrow A$. In Section~$3$ we
 compare the Lie algebra homomorphisms $\theta_{\nu}$ for different presentations
$\nu$ of $A$, using the quiver $\Gamma$. Finally, in Section~$4$ we
 prove Theorem~\ref{thm1}.

This text is part of the author's thesis (\cite{lemeur_thesis}) made
at Universitï¿½~Montpellier~II under the supervision of Claude Cibils.

\section{Preliminaries}
\subsection{Terminology and notations for quivers} 
Let $Q$ be a quiver. We write 
$Q_0$ and $Q_1$ for the set of vertices and of arrows,
respectively. We read (oriented) paths from the right to the left,
that is, we view a path $u$ as a morphism and the concatenation $vu$
of two paths $u$ and $v$ such that the source of $v$ equals the target
of $u$ as a
composition of morphisms. Given $x\in Q_0$, the \emph{trivial} path (of
length $0$, with source and target equal to $x$) is denoted by $e_x$.
 Two
paths are called \emph{parallel} if they have the same source and
the same target.
 An
\emph{oriented cycle} in $Q$ is a non trivial path whose source
and target are equal. If $\alpha\in Q_1$ we consider its formal
inverse $\alpha^{-1}$ with source and target equal to the target and
the source of $\alpha$, respectively. Hence, we get the double quiver
$\overline{Q}$ such that $\overline{Q}_0=Q_0$ and
$\overline{Q}_1=Q_1\cup\{\alpha^{-1}\ |\ \alpha\in Q_1\}$. Then, a
\emph{walk} in $Q$ is exactly an oriented path in $\overline{Q}$.
Given a walk
$\gamma=\alpha_n^{\varepsilon_n}\ldots\alpha_1^{\varepsilon_1}$ (with
$\alpha_i\in Q_1$, $\varepsilon_i\in\{\pm 1\}$), its inverse
$\gamma^{-1}$ is by definition $\alpha_1^{-\varepsilon_1}\ldots\alpha_n^{-\varepsilon_n}$.

\subsection{Presentations  by quiver and admissible relations} Let $Q$ be a  quiver.
Its \emph{path algebra} $kQ$ is
the $k$-algebra whose basis as a $k$-vector space is the set of paths
in $Q$ (including the trivial paths), and whose product is bilinearly
induced by the concatenation of paths (if $u,v$ are two paths such
that the source of $v$ is different from the target of $u$, then we
set $vu=0$). The unit of $kQ$ is
$\sum\limits_{x\in Q_0}e_x$ and $kQ$ is finite dimensional if and only if
$Q$ is finite (that is $Q_0$ and $Q_1$ are finite) and has no
oriented cycles. We let $kQ^+$ be the ideal of $kQ$ generated by the
arrows.

An \emph{admissible ideal} of $kQ$ is an
ideal $I$ such that $(kQ^+)^N\subseteq I\subseteq (kQ^+)^2$ for some
$N\geqslant 2$.
In such a case, the elements of $I$ are called \emph{relations} and,
following \cite{martinezvilla_delapena}, a \emph{minimal relation}
of $I$ is a relation
$\sum\limits_{i=1}^st_iu_i\neq 0$ such that
 $t_1,\ldots,t_s\in k^*$ and $u_1,\ldots,u_s$ are pairwise
  distinct paths in $Q$,
and such that there is no non empty proper subset $S\subset
\{1,\ldots,s\}$ satisfying $\sum\limits_{i\in
  S}t_iu_i\in I$.
In such a case,  $u_1,\ldots,u_s$ are necessarily parallel.
Note that $I$ is generated by its minimal relations.

Recall (see \cite{ars}) that any \fd \kalg $A$ is Morita equivalent to a
basic one.
If $A$ is
basic, then  there exists a unique quiver $Q$, the \emph{ordinary
  quiver} of $A$, and a
surjective \kalg homomorphism $\nu\colon kQ\twoheadrightarrow A$ whose
kernel is an admissible ideal of $kQ$.  Also, $\{\nu(e_x)\ |\
x\in Q_0\}$ is a complete set of primitive orthogonal idempotents of $A$. The
homomorphism $\nu$ is called a \emph{ presentation (by quiver and
  admissible relations)}. We have $A\simeq kQ/\ker(\nu)$ and $A$ is connected if and
only if $Q$ is connected.\\

\subsection{Presentation of $\hh^1(A)$} Let $A$ be a \bfda and let
$\{e_1,\ldots,e_n\}$ be a complete set of primitive orthogonal
idempotents. A \emph{unitary
  derivation} (\cite{delapena_saorin}) is a $k$-linear map $d\colon A\to A$ such
that $d(ab)=ad(b)+d(a)b$ for any $a,b\in A$ and such that $d(e_i)=0$
for every $i$. Let $Der_0(A)$ be set of unitary
derivations. It is a Lie algebra for the commutator. In the sequel, all derivations will be unitary. So
we shall call them derivations. Let
$E:=\{\sum\limits_{i=1}^nt_ie_i\ |\ t_1,\ldots,t_n\in k\}$. Then $E$ is a semi-simple
subalgebra of $A$ and $A=E\oplus\mathfrak{r}$ where $\mathfrak{r}$ is
the radical of $A$. Let $Int_0(A):=\{\delta_e\colon A\to A,\ a\in
A\mapsto ea-ae\ |\ e\in E\}$, this is an ideal of
$Der_0(A)$. Throughout this text, we shall use the following
presentation proved in \cite{delapena_saorin}:
\begin{thm}
\label{thm1.6}
(\cite{delapena_saorin}) $\hh^1(A)\simeq Der_0(A)/ Int_0(A)$
as Lie algebras.
\end{thm}
In the following lemma, we collect some useful properties on
derivations.
\begin{lem}
\label{lem1.1}
Let $d\in Der_0(A)$, then $d(e_jAe_i)\subseteq e_jAe_i$.  Assume that the
ordinary quiver $Q$ of $A$ has no oriented cycles, then $d(\mathfrak{r})\subseteq
\mathfrak{r}$ and $d(\mathfrak{r}^2)\subseteq \mathfrak{r}^2$.
\end{lem}
\noindent{\textbf{Proof:}} Since $d$ is unitary and since $Q$ has no
oriented cycles, we have $d(\mathfrak{r})\subseteq\mathfrak{r}$. So, $d(\mathfrak{r}^2)\subseteq\mathfrak{r}^2$.
\hfill$\square$\\

If $\psi\colon A\xrightarrow{\sim} A$ is a \kalg automorphism such
that $\psi(e_i)=e_i$ for every $i$, then the map $d\mapsto
\psi\circ d\circ \psi^{-1}$ induces a Lie algebra automorphism of
$\hh^1(A)$, denoted by $\psi_*\colon \hh^1(A)\xrightarrow{\sim} \hh^1(A)$.\\

\subsection{Fundamental groups of presentations} Let $(Q,I)$ be a
\emph{bound quiver} (that is, $Q$ is a finite quiver and $I$ is an
admissible ideal of $kQ$). The \emph{homotopy relation} $\sim_I$ was
defined in \cite{martinezvilla_delapena} as the equivalence class
on the set of walks in $Q$ generated by the following properties:
\begin{enumerate}
\item[($1$)] $\alpha\alpha^{-1}\sim_I e_y$ and
  $\alpha^{-1}\alpha\sim_I e_x$ for any arrow $\alpha$ with source $x$
  and target $y$,
\item[($2$)] $wvu\sim_Iwv'u$ if $w,v,v',u$ are walks such that the
  concatenations $wvu$
  and $wv'u$ are well defined and such that $v\sim_Iv'$,
\item[($3$)]$u\sim_I v$ if $u$ and $v$ are paths in a same
  minimal relation (with a non zero scalar). 
\end{enumerate}
Note that if $r_1,\ldots,r_t$ are minimal relations generating
the ideal $I$, then the condition ($3$) above may be
replaced by the following one (\cite{bardzell_marcos}):
\begin{enumerate}
\item[($3'$)] $u\sim_I v$ if $u$ and $v$ are paths in $Q$ appearing in
  $r_i$ (with a non zero scalar) for some $i\in\{1,\ldots,t\}$.
\end{enumerate}
The $\sim_I$-equivalence class of a walk $\gamma$ is be denoted  by
$[\gamma]_I$. Let $x_0\in Q_0$, following
\cite{martinezvilla_delapena}, the set of $\sim_I$-equivalence classes
of walks with source and target $x_0$ is denoted by
$\pi_1(Q,I,x_0)$. The concatenation of walks endows this set with a
group structure whose unit is $[e_{x_0}]_I$. This group is called the
\emph{fundamental group} of $(Q,I)$ at $x_0$. If $Q$ is connected,
then the isomorphism class of $\pi_1(Q,I,x_0)$ does not depend on the
choice $x_0$. In such a case, we write $\pi_1(Q,I)$ for
$\pi_1(Q,I,x_0)$. If $A$ is a \bcfda and if $\nu\colon
kQ\twoheadrightarrow A$ is a presentation, the group
$\pi_1(Q,\ker(\nu))$  is called \emph{the fundamental group of
  the presentation $\nu$}. The following example shows that two
presentations of $A$ may have non isomorphic fundamental groups.
\begin{ex}
\label{ex1.7}
Let $A=kQ/I$ where $Q$ is the quiver:
\xymatrix{
1 \ar@/^/[r]^b \ar@/_/[r]_a & 2 \ar[r]^c& 3}
and $I=<ca>$. Set $x_0=1$. Then $\pi_1(Q,I)\simeq\mathbb{Z}$ is
generated by $[b^{-1}a]_I$. On the other hand, $A\simeq kQ/J$
where $J=<ca-cb>$, and $\pi_1(Q,J)$ is the trivial group.
\end{ex}

In the sequel we shall use the following technical lemma.
\begin{lem}
\label{lem1.2}
Let $(Q,I)$ be a bound quiver where $Q$ has no oriented cycles and let
$d\colon kQ\to kQ$ be a linear map such that
$d(I)\subseteq I$, and $d(u)=t_uu$ for some $t_u\in k$, for any
path $u$. Let $\equiv_I$ be the equivalence relation on the set of paths in $Q$
generated by the condition ($3$) defining $\sim_I$.
Then, the following implication holds for any paths $u,v$:
\begin{equation}
u\equiv_Iv\ \text{implies}\ t_u=t_v\ .\notag
\end{equation}
\end{lem}
\noindent{\textbf{Proof:}} We use a non multiplicative version
of Grï¿½bner bases (\cite{adams_loustaunau}, see also
\cite{farkas_feustel_green}). Fix an arbitrary total order
$u_1<\ldots<u_N$ on the set of paths in $Q$ and let
$(u_1^*,\ldots,u_N^*)$ be the basis of $\hom_k(kQ,k)$ dual to
$(u_1,\ldots,u_N)$. Following \cite[Sect. 1]{lemeur2}, the Grï¿½bner
basis of $I$ is the unique basis $(r_1,\ldots,r_t)$ defined by the
three following properties:
\begin{enumerate}
\item[($i$)] $r_j\in u_{i_j}\  +\ Span(u_i\ ;\ i<i_j)$ for some $i_j$, for
  every $j$,
\item[($ii$)] $u_{i_j}^*(r_{j'})=0$ unless $j=j'$,
\item[($iii$)] $i_1<\ldots<i_t$.
\end{enumerate}
It follows from these properties that:
\begin{enumerate}
\item[$(iv$)]
$r=\sum\limits_{j=1}^tu_{i_j}^*(r)r_j$ for any $r\in I$.
\end{enumerate}
Recall from \cite[Sect. 1]{lemeur2} that $r_1,\ldots,r_t$ are minimal relations of
$I$ so that $\equiv_I$ is generated by the property ($3'$) defining
$\sim_I$. So we only need to prove that
that $d(r_j)\in k.r_j$ for any $j$. We proceed by
induction on $j\in\{1,\ldots,t\}$. By assumption on $d$ and thanks to
($i$), we have
$d(r_1)\in I\cap Span(u_i\ ;\ i\leqslant i_1)$. Hence, ($iii$) and
($iv$) imply
that $d(r_1)\in
k.r_1$. Let $j\in\{1,\ldots,d-1\}$ and assume that $d(r_1)\in
k.r_1,\ldots,d(r_j)\in k.r_j$. By assumption on $d$ and thanks to
($i$) and ($ii$), we have
$d(r_{j+1})\in I\cap Span(u_i\ ;\ i\leqslant i_{j+1})$ and
$u_{i_l}^*(d(r_{j+1}))=0$ if $l\leqslant j$. So, ($iii$) and ($iv$)
imply that $d(r_{j+1})\in k.r_{j+1}$. This finishes the induction and
proves the lemma.\hfill$\square$\\ 

\subsection{Comparison of fundamental groups} Let $A$ be a \bcfda with
ordinary quiver $Q$ without oriented cycles. We defined the
transvections in the introduction. A \emph{dilatation} (\cite{lemeur}) is an
automorphism $D\colon kQ\xrightarrow{\sim} kQ$ such that $D(e_i)=e_i$
for any $i$ and such that $D(\alpha)\in k.\alpha$ for any $\alpha\in
Q_1$. The following proposition will be useful in the sequel, it was proved in \cite{lemeur2}:
\begin{prop}
\label{prop1.3}(\cite[Prop. 2.5]{lemeur2})
Let $I$ be an admissible
ideal of $kQ$, let $\varphi$ be an automorphism of $kQ$ and let
$J=\varphi(I)$. If $\varphi$ is a dilatation,
then $\sim_I$ and $\sim_J$ coincide. Assume that
$\varphi=\varphi_{\alpha,u,\tau}$:
\begin{enumerate}
\item[-] If $\alpha\sim_Iu$ and $\alpha\sim_Ju$ then $\sim_I$ and
$\sim_J$ coincide.
\item[-] If $\alpha\not\sim_Iu$ and $\alpha\sim_J u$ then $\sim_J$
is generated by $\sim_I$ and  $\alpha\sim_Ju$.
\item[-] If $\alpha\not\sim_Iu$ and $\alpha\not\sim_Ju$ then $I=J$
and $\sim_I$ and $\sim_J$ coincide.
\end{enumerate}
In particular, if $\alpha\sim_J u$, then the identity map on the set
of walks in $Q$ induces a surjective group homomorphism
$\pi_1(Q,I)\twoheadrightarrow \pi_1(Q,J)$.
\end{prop}
Here \textit{generated} means: generated as an equivalence relation on
the set of walks in $Q$, and satisfying the conditions ($1$) and ($2$)
 in the definition of the homotopy relation. 
If $I,J$ are admissible ideals such that there exists
$\varphi_{\alpha,u,\tau}$ satisfying $J=\varphi_{\alpha,u,\tau}(I)$,
$\alpha\not\sim_I u$ and $\alpha\sim_Ju$, then we say that
\emph{$\sim_J$ is a direct successor of $\sim_I$}.
Proposition~\ref{prop1.3} allows one to define a
quiver $\Gamma$ associated to $A$ as follows (\cite[Def. 4.1]{lemeur}):
\begin{enumerate}
\item[-] $\Gamma_0=\{\sim_I\ |\ \text{$I$ is an admissible ideal of
$kQ$ such that $A\simeq kQ/I$}\}$,
\item[-] there is an arrow $\sim\to\sim'$ if $\sim_J$ is a direct
successor of $\sim_I$.
\end{enumerate}

\begin{ex}
\label{ex1.8}
Let $A$ be as in Example~\ref{ex1.7}, then $J=\varphi_{\alpha,cb,1}(I)$ and
$\Gamma$ is equal to $\sim_I\to \sim_J$.
\end{ex}

The quiver $\Gamma$ is finite, connected and has not oriented cycles
(\cite[Rem. 3, Prop. 4.2]{lemeur}). Moreover, if $\Gamma$ has a unique
source $\sim_{I_0}$ (that is, a vertex with no arrow ending at it) then
the fundamental group of any admissible presentation of $A$ is a
quotient of $\pi_1(Q,I_0)$. It was proved in \cite{lemeur2} and
\cite{lemeur3} that $\Gamma$ has a unique source under one of the
 hypotheses ($H_1$) or ($H_2$) presented in the introduction.
Moreover, the hypotheses ($H_1$) and ($H_2$) both ensure the following proposition which
will be particularly useful to prove Theorem~\ref{thm1}.
\begin{prop}
\label{prop1.4}
(\cite[Lem. 4.3]{lemeur2} and \cite[Prop. 4.3]{lemeur3})
Assume that at least one of the two hypotheses ($H_1$) or ($H_2$) is satisfied.
Let $\sim_{I_0},\sim_I\in\Gamma$, where $\sim_{I_0}$ is the unique
source of $\Gamma$. Then there exist a dilatation $D$ and a sequence of
transvections
$\varphi_{\alpha_1,u_1,\tau_1},\ldots,\varphi_{\alpha_l,u_l,\tau_l}$
such that:
\begin{enumerate}
\item[-]
$I=D\varphi_{\alpha_l,u_l,\tau_l}\ldots\varphi_{\alpha_1,u_1,\tau_1}(I_0)$.
\item[-] If we set
  $I_i:=\varphi_{\alpha_i,u_i,\tau_i}\ldots\varphi_{\alpha_1,u_1,\tau_1}(I_0)$,
  then $\alpha_i\sim_{I_i}u_i$ for every $i$.
\end{enumerate}
If $\sim_I=\sim_{I_0}$, then
$\sim_{I_0},\sim_{I_1},\ldots,\sim_{I_l},\sim_I$ coincide.
\end{prop}

\subsection{Comparison of the fundamental groups and the Hochschild
  cohomology} Let $A$ be a \bcfda. Assume that the ordinary quiver
  $Q$ of $A$ has no oriented cycles. Let $x_0\in Q_0$ and fix a maximal
  tree $T$ of $Q$,
  that is, a subquiver of $Q$ such that $T_0=Q_0$ and such that the
  underlying graph of $T$ is a tree. With these data, Assem and
  de~la~Pe\~na have defined  an
  injective homomorphism of abelian
  groups $\theta_{\nu}\colon \hom(\pi_1(Q,\ker(\nu)),k^+)\hookrightarrow
  \hh^1(A)$ associated to any admissible presentation $\nu\colon
  kQ\twoheadrightarrow A$ (\cite{assem_delapena}). We recall the definition of
  $\theta_{\nu}$ and refer the reader to \cite{assem_delapena} for
  more details. For any $x\in Q_0$ there exists a unique walk
  $\gamma_x$ in $T$ with source $x_0$, with target $x$ and of minimal length
 for these properties. Let $\nu\colon kQ\twoheadrightarrow A$
  be an admissible presentation and let $f\in \hom(\pi_1(Q,\ker(\nu)),k^+)$ be a
  group homomorphism. Then, $f$ defines a derivation $\widetilde{f}\colon
  A\to A$ as follows:
  $\widetilde{f}(\nu(u))=f([\gamma_y^{-1}u\gamma_x]_{\sim_{\ker(\nu)}})\
  \nu(u)$ for any path $u$ with source $x$ and target $y$. The
  following proposition was proved in \cite{assem_delapena}:
\begin{prop}
(\cite{assem_delapena})
The map $f\mapsto \widetilde{f}$ induces an injective map of
abelian groups:
\begin{equation}
\theta_{\nu}\colon \hom(\pi_1(Q,\ker(\nu)),
k^+)\hookrightarrow \hh^1(A)\ .\notag
\end{equation}
\end{prop}
Note that $\theta_{\nu}$ is not surjective in
general. Indeed, if $A$ is the path algebra of the Kronecker quiver, then
$\ker(\nu)=0$, $\dim_k\
\im(\theta_{\nu})=1$, and
$\dim_k\ \hh^1(A)=3$. Note also that despite its definition, the
homomorphism $\theta_{\nu}$ does not depend on the choice of
$T$. Indeed, let $T'$ be another maximal tree, thus defining the walk
$\gamma_x'$ of minimal length in $T'$ with source $x_0$ and target $x$, for every vertex
$x$. Given a group homormorphism $f\colon \pi_1(Q,\ker(\nu))\to
k^+$ there is a new derivation $\hat{f}\colon A \to A$ (instead of
$\widetilde{f}$) obtained by applying the previous
construction to $T'$ (instead of to $T$), that is
$\hat{f}(\nu(u))=f([\gamma_y'^{-1}u\gamma'_x]_{\ker(\nu)})\nu(u)$ for every path $u$ in
$Q$ from $x$ to $y$. Now let $e=\sum\limits_{x\in
  Q_0}f([\gamma_x'^{-1}\gamma_x]_{\ker(\nu)})\,e_x\in A$. It is easily
checked that $\hat{f}-\widetilde{f}$ is the inner derivation
associated to $e$. In particular, $\widetilde{f}$ and $\hat{f}$ have
equal images in $\hh^1(A)$. So the construction of $\theta_{\nu}$ does
not depend on the choice of the maximal tree $T$.
 
The product in $k$ endows $\hom(\pi_1(Q,\ker(\nu)),k^+)$ with a
commutative $k$-algebra structure. So it is also an abelian Lie
algebra for the commutator. The following lemma proves that
$\theta_{\nu}$ preserves this
structure. The proof is just a direct computation, so we omit it.
\begin{lem}
$\theta_{\nu}\colon \hom(\pi_1(Q,\ker(\nu)),
k^+)\hookrightarrow \hh^1(A)$ is a Lie algebra homomorphism. In particular,
$\im(\theta_{\nu})$ is an abelian Lie subalgebra of $\hh^1(A)$.
\end{lem}

Throughout this text, $A$ will be a basic connected finite dimensional
$k$-algebra with ordinary quiver $Q$ without oriented cycles
($Q_0=\{1,\ldots,n\}$). We fix
a
complete set $\{e_1,\ldots,e_n\}$  of primitive
orthogonal idempotents of $A$. 
So $A=E\oplus\mathfrak{r}$, where $E=k.e_1\oplus\ldots\oplus k.e_n$
and $\mathfrak{r}$ is the radical of $A$. Without
loss of generality, we assume that any presentation $\nu\colon
kQ\twoheadrightarrow A$ is such that $\nu(e_i)=e_i$. Finally, in order to
use the Lie algebra homomorphisms $\theta_{\nu}$, we fix a maximal tree $T$ in $Q$.

\section{Diagonalizability in $\hh^1(A)$}
The aim of this section is to prove some useful properties on the
subspaces $\im(\theta_{\nu})$ in terms of diagonalizability in
$\hh^1(A)$. Note that diagonalizability was introduced for
derivations of $A$ in \cite{farkas_green_marcos}.
For short, a \emph{basis} of $A$ is a basis $\c B$ of the $k$-vector
space $A$ such that:
$\c B\subseteq \bigcup\limits_{i,j}e_jAe_i$,
such that $\{e_1,\ldots,e_n\}\subseteq \c B$,
and such that $\c B\backslash\{e_1,\ldots,e_n\}\subseteq\mathfrak{r}$.
 Note the following link between bases
and presentations of $A$:
\begin{enumerate}
\item[-] If $\nu\colon kQ\twoheadrightarrow A$ is a
  presentation of $A$, then there exists a basis $\c B$ such that
  $\nu(\alpha)\in\c B$ for any $\alpha\in Q_1$ and such that any element
  of $\c B$ is of the form $\nu(u)$ with $u$ a path in $Q$. We  say
  that this basis $\c B$ is \emph{adapted} to $\nu$.
\item[-] If $\c B$ is a basis of $A$, then there exists a
  presentation $\nu\colon kQ\twoheadrightarrow A$ such that
  $\nu(\alpha)\in\c B$ for any $\alpha\in Q_1$. We  say that the
  presentation $\nu$ is adapted to $\c B$.
\end{enumerate}
The property of being diagonalizable (as a linear map) is stable under the sum with an inner
derivation as the following lemma shows. The proof is immediate.
\begin{lem}
\label{lem2.1}
Let $u\colon A\to A$ be a linear map, let $e\in E$ and let $\c B$ be a basis of
$A$. Then $u$ is diagonal with respect to the basis $\c B$ if and only if the
same holds for $u+\delta_e$.
\end{lem}

The preceding lemma justifies the following definition.
\begin{definition}
\label{def2.2}
Let $f\in \hh^1(A)$ and let $d$ be a derivation representing $f$. Then
$f$ is called \emph{diagonalizable} (and \emph{diagonal} with respect to a basis $\c B$ of $A$) if and only if $d$ is diagonalizable (and diagonal with respect to $\c B$, respectively).

The subset $D\subseteq \hh^1(A)$ is called diagonalizable if and only
if any there exists a basis $\c B$ of $A$ such that any $f\in D$ is diagonal with respect to $\c B$.
\end{definition}

The following proposition gives a criterion for a subset $D\subset
\hh^1(A)$ to be diagonalizable.
\begin{prop}
\label{prop2.3}
Let $D\subseteq \hh^1(A)$. Then,
$D$ is diagonalizable if and only if
every element of $D$ is diagonalizable and $[f,f']=0$ for any $f,f'\in
D$.
\end{prop}
\noindent{\textbf{Proof:}} Clearly, if $D$ is diagonalizable, then so is every element of $D$ and $[f,f']=0$ for every $f,f'\in D$. We prove the converse.
For each $f\in D$, let $d_f$ be a derivation
representing $f$. So $d_f$ is diagonal with respect to some basis and it suffices to prove that this basis may be assumed to be the same for all $f\in D$.
Note that $d_f$ induces a diagonalizable linear map $d_f\colon
e_j\mathfrak{r}e_i\to e_j\mathfrak{r}e_i$, for every $i,j$ (see
Lemma~\ref{lem1.1}). Also, for every $f,f'\in D$, there exist scalars
$t_i^{(f,f')}\in k$, for $i\in\{1,\ldots,n\}$, such that
$[d_f,d_{f'}]$ is the inner derivation
$\delta_{e^{(f,f')}}$, where
$e^{(f,f')}=\sum\limits_{i=1}^nt_i^{(f,f')}e_i$. Now, let
$i,j\in\{1,\ldots,n\}$. Then, given $f,f'\in D$, we have two
diagonalizable maps $d_f,d_{f'}\colon e_j\mathfrak{r}e_i\to
e_j\mathfrak{r}e_i$ whose commutator is equal to
$(t_j^{(f,f')}-t_i^{(f,f')})Id_{e_j\mathfrak{r}e_i}$. So this
commutator must be zero. This shows that there exists a basis $\c
B_{i,j}$ of $e_j\mathfrak{r}e_i$ for which $d_f\colon
e_j\mathfrak{r}e_i\to e_j\mathfrak{r}e_i$ has a diagonal matrix. So any $f\in D$ is diagonal with respect to the basis
$\c B=\{e_1,\ldots,e_n\}\cup\bigcup\limits_{i,j}\c B_{i,j}$ which does not depend on $f$. This proves that
$D$ is diagonalizable.\hfill$\square$\\

Our main example of diagonalizable subspace of
$\hh^1(A)$ is $\im(\theta_{\nu})$:
\begin{prop}
\label{prop2.4}
Let $\nu\colon kQ\twoheadrightarrow A$ be a
presentation. Then, $\im(\theta_{\nu})$ is diagonalizable.
\end{prop}
\noindent{\textbf{Proof:}} Let $\c B$ be a basis of $A$ adapted to $\nu$
and let $I=\ker(\nu)$. Then $\theta_{\nu}(f)$ is diagonal with respect to $\c B$, for every $f\in \hom(\pi_1(Q,I),k^+)$.\hfill$\square$\\ 

In this section, we aim at proving that any 
diagonalizable subset of $\hh^1(A)$ is contained in $\im(\theta_{\nu})$
for some presentation $\nu$.
It was proved in \cite{farkas_green_marcos} that any diagonalizable
derivation (with
suitable technical conditions) defines an element of $\hh^1(A)$ lying
in $\im(\theta_{\nu})$ for some $\nu$. We will use
the following similar result.
\begin{lem}
\label{lem2.5}
Let $f\in \hh^1(A)$ be diagonalizable. Let $\c B$ be a basis with respect to which $f$ is diagonal. Let
$\nu\colon kQ\twoheadrightarrow A$ be a presentation
adapted to $\c B$. Then $f\in \im(\theta_{\nu})$.
\end{lem}
\noindent{\textbf{Proof:}} Let $I=\ker(\nu)$ and let $d\colon A\to A$
be a derivation representing $f$. We  set $\overline{r}:=\nu(r)$, for any
$r\in kQ$. Let $\alpha\in Q_1$. By assumption on $\c
B$, there exists $t_{\alpha}\in k$
such that $d(\overline{\alpha})=t_{\alpha}\overline{\alpha}$. Let
$t_u:=t_{\alpha_1}+\ldots+t_{\alpha_n}$, for any path
$u=\alpha_n\ldots\alpha_1$ (with $\alpha_i\in Q_1$).  So  $d(\overline{u})=t_u\overline{u}$,
because $d$ is a derivation. More generally, if
$\gamma=\alpha_n^{\e_n}\ldots\alpha_1^{\e_1}$ is a walk in $Q$ (with
$\alpha_i\in Q_1$), let us set
$t_{\gamma}:=\sum\limits_{i=1}^n(-1)^{\e_i}t_{\alpha_i}$, with the
convention that $t_{\gamma}=0$ if $\gamma$
is trivial. We now to prove that the map $\gamma\mapsto
t_{\gamma}$ defines a group homomorphism $g\colon \pi_1(Q,I)\to
k^+,[\gamma]_I\mapsto t_{\gamma}$ and
that $f=\theta_{\nu}(g)$.

First, we prove that the group homomorphism $g\colon \pi_1(Q,I)\to k^+$ is
well defined. By
definition of the scalar $t_{\gamma}$, we have:
\begin{enumerate}
\item[($i$)] $t_{e_{x}}=0$ for any $x\in Q_0$ and
$t_{\gamma'\gamma}=t_{\gamma'}+t_{\gamma}$ for any walks
$\gamma,\gamma'$ such that the walk $\gamma'\gamma$ is defined.
\item[($ii$)] $t_{\alpha^{-1}\alpha}=t_{e_x}$ and
$t_{\alpha\alpha^{-1}}=t_{e_y}$ for any arrow
$x\xrightarrow{\alpha}y\in Q_1$.
\item[($iii$)] $t_{wvu}=t_{wv'u}$ for any walks $w,v,v',u$ such that 
$t_v=t_{v'}$, and such that the walks $wvu,\,wv'u$ are defined.
\end{enumerate}
In order to prove that $g$ is well defined, it only remains to prove
that $t_u=t_v$ whenever $u,v$ are paths in $Q$ appearing in
the same minimal relation of $I$ (with non zero scalars). For this purpose, let $d'\colon kQ\to kQ$ be
the linear map such that $d'(u)=t_uu$ for any path $u$ in $k$. Thus, $d\circ \nu=\nu\circ
d'$. In particular, $d'(I)\subseteq I$. So we may apply
Lemma~\ref{lem1.2} to $d'$ and deduce that:
\begin{enumerate}
\item[($iv$)] $t_u=t_v$ if $u,v$ are
paths in $Q$ lying in the support of a same minimal relation of $I$.
\end{enumerate}
From ($ii$), ($iii$) and ($iv$) we deduce that we have a well defined map
$g\colon \pi_1(Q,I)\to k,[\gamma]_I\mapsto t_{\gamma}$. Moreover, ($i$) proves that $g$ is a group homomorphism.

Now we prove that $f=\theta_{\nu}(g)$. For any path $u$ with
source $x$ and target $y$, we have
$g([\gamma_y^{-1}u\gamma_x]_I)=t_u-t_{\gamma_y}+t_{\gamma_x}$. Hence,
$\theta_{\nu}(g)\in \hh^1(A)$ is represented by the derivation
$\widetilde{g}\colon A\to A$ such that
$\widetilde{g}(\overline{u})=(t_u-t_{\gamma_y}+t_{\gamma_x})\overline{u}$ for any path $u$ with
source $x$ and target $y$. Let us set $e:=\sum\limits_{x\in
Q_0}t_{\gamma_x}e_x\in E$. Therefore, $\widetilde{g} +\delta_e=d$. This proves that
$f=\theta_{\nu}(g)$.\hfill$\square$\\

Now we can state the main result of this section. It is a direct
consequence of Proposition~\ref{prop2.4} and of Lemma~\ref{lem2.5}. 
\begin{prop}
\label{prop2.6}
Let $D\subseteq \hh^1(A)$. Then $D$ is diagonalizable if and only if
there exists a
presentation $\nu\colon kQ\twoheadrightarrow A$ such that
$D\subseteq \im(\theta_{\nu})$.
\end{prop}

Remark that Lemma~\ref{lem2.5} also gives a sufficient condition for
$\theta_{\nu}$ to be an isomorphism. Recall that $A$ is called
constricted if and only if $\dim\ e_yAe_x=1$ for any arrow $x\to y$
(this implies that $Q$ has no multiple arrows). In
\cite{bardzell_marcos} it was proved that for such an algebra, two
different presentations have the same fundamental group.
\begin{prop}
\label{prop2.7}
Assume that $A$ is constricted. Let $\nu\colon kQ\twoheadrightarrow A$
be any presentation of $A$. Then $\theta_{\nu}\colon
\hom(\pi_1(Q,I),k^+)\to \hh^1(A)$ is an isomorphism. In particular,
$\hh^1(A)$ is an abelian Lie algebra.
\end{prop}
\noindent{\textbf{Proof:}} Since $\theta_{\nu}$ is one-to-one, we only
need to prove that it is onto. Let $\c B$ be a basis of $A$ adapted to
$\nu$, let $f\in \hh^1(A)$ and let $d\colon A\to A$ be a derivation
representing $f$. Let $x\xrightarrow{\alpha}y$ be an arrow. Then
$e_yAe_x=k.\nu(\alpha)$ so that there exists $t_{\alpha}\in k$ such
that $d(\nu(\alpha))=t_{\alpha}\nu(\alpha)$. Let $u=\alpha_n\ldots\alpha_1$ be any path in
$Q$ (with $\alpha_i\in Q_1$). Since $d$ is a derivation, we have
$d(\nu(u))=(t_{\alpha_1}+\ldots+t_{\alpha_n})\nu(u)$. As a
consequence, $d$ is diagonal with respect to $\c B$. Moreover, $\nu$ is
adapted to $\c B$. So Lemma~\ref{lem2.5}
proves that $f\in \im(\theta_{\nu})$. This proves that $\theta_{\nu}$
is an isomorphism. So $\hh^1(A)$ is abelian.\hfill$\square$\\

\section{Comparison of $\im(\theta_{\nu})$ and $\im(\theta_{\mu})$ for
  different presentations $\mu$ and $\nu$ of $A$}
If two
presentations $\nu$ and $\mu$ of $A$ are related by a transvection or a
dilatation, then there is a simple relation between the associated
fundamental groups (see Proposition~\ref{prop1.3}). In this section, we  
compare $\theta_{\nu}$ and $\theta_{\mu}$.
We first compare $\theta_{\nu}$ and $\theta_{\mu}$ when
$\mu=\nu\circ D$ with $D$ a dilatation. Recall that if $J=D(I)$ with $D$ a dilatation, then
$\sim_I$ and $\sim_J$ coincide, so that $\pi_1(Q,I)=\pi_1(Q,J)$.
\begin{prop}
\label{prop3.1}
Let $\nu\colon kQ\twoheadrightarrow A$ be a presentation, let $D\colon
kQ\xrightarrow{\sim}kQ$ be a dilatation.
Let $\mu:=\nu\circ D\colon kQ\twoheadrightarrow A$. Let $I=\ker(\mu)$
and $J=\ker(\nu)$, so that $J=D(I)$. Then $\theta_{\mu}=\theta_{\nu}$.
\end{prop}
\noindent{\textbf{Proof:}} Let $f\in \hom(\pi_1(Q,I),k^+)$. Then,
$\theta_{\nu}(f)$ and $\theta_{\mu}(f)$ are represented by the
derivations $d_1$ and $d_2$ respectively, such that for any arrow $x\xrightarrow{\alpha}y$:
\begin{align}
&d_1(\nu(\alpha))=f([\gamma_y^{-1}\alpha\gamma_x]_J)\ \nu(\alpha)\notag\\
&d_2(\mu(\alpha))=f([\gamma_y^{-1}\alpha\gamma_x]_I)\ \mu(\alpha)\ .\notag
\end{align}
Therefore, $d_1(\nu(\alpha))=d_2(\mu(\alpha))$ because $D$ is a
dilatation and because $\sim_I$ and $\sim_J$ coincide. This implies
that $d_1=d_2$ and $\theta_{\nu}(f)=\theta_{\mu}(f)$.\hfill$\square$\\

The following example shows that
Proposition~\ref{prop3.1} does not necessarily hold true if $\nu$ and $\mu$ are two
presentations of $A$ such that $\sim_{\ker(\nu)}$ and $\sim_{\ker(\mu)}$
coincide.
\begin{ex}
\label{ex3.4}
Assume that $\carac(k)=2$ and let $A=kQ/I$ where $Q$ is the quiver:
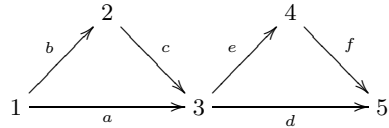
\begin{equation}
\xymatrix{
& 2 \ar@{->}[rd]^c & & 4\ar@{->}[rd]^f&\\
1\ar@{->}[ru]^b \ar@{->}[rr]_a & & 3 \ar@{->}[ru]^e \ar@{->}[rr]_d & &5
}\notag
\end{equation}
and $I=<da,fecb,fea+dcb>$. Let $T$ be the maximal tree such that
$T_1=\{b,c,e,f\}$. Let $\nu\colon kQ\twoheadrightarrow A=kQ/I$ be
the natural projection. Let $\psi:=\varphi_{a,cb,1}\varphi_{d,fe,1}$.
Thus, $I=\psi(I)$. Let $\mu:=\nu\circ \psi\colon
kQ\twoheadrightarrow A$ so that $\ker(\mu)=\ker(\nu)=I$. Observe that
$\pi_1(Q,I)$ is the infinite cyclic group with generator
$[b^{-1}c^{-1}a]_I$. So let $f\colon \pi_1(Q,I)\to k^+$ be the unique group
homomorphism such that $f([b^{-1}c^{-1}a]_I)=1$. Then $\theta_{\nu}(f)$ is
represented by the following derivation:
\begin{equation}
\begin{array}{rrclc}
d_1 \colon & A & \longrightarrow & A\\
& \nu(x) & \longmapsto & \nu(x)\ \text{if $x\in \{a,d\}$}\\
& \nu(x) & \longmapsto & 0\ \text{if $x\in \{b,c,e,f\}$}&.
\end{array}\notag
\end{equation}
On the other hand, $\theta_{\mu}(f)$ is represented by the derivation:
\begin{equation}
\begin{array}{rrclc}
d_2 \colon & A & \longrightarrow & A\\
& \nu(a) & \longmapsto & \nu(a)+\nu(cb)\\
& \nu(d) & \longmapsto & \nu(d)+\nu(fe)\\
& \nu(x) & \longmapsto & 0\ \text{if $x\in\{b,c,e,f\}$}&.
\end{array}\notag
\end{equation}
It is easy to verify that $d_2-d_1$ is not an inner derivation. Hence,
$\theta_{\nu}\neq\theta_{\mu}$.
\end{ex}

Now we compare $\theta_{\nu}$ and $\theta_{\mu}$ when
$\mu=\nu\circ\varphi_{\alpha,u,\tau}$ and when the identity map on the
set of walks in $Q$ induces a surjective group homomorphism
$\pi_1(Q,\ker(\nu))\twoheadrightarrow \pi_1(Q,\ker(\mu))$.
\begin{prop}
\label{prop3.2}
Let $\nu\colon kQ\twoheadrightarrow A$ be a presentation, let
$\varphi_{\alpha,u,\tau}\colon kQ\xrightarrow{\sim} kQ$ be a
transvection and let $\mu:=\nu\circ\varphi_{\alpha,u,\tau}\colon
kQ\twoheadrightarrow A$. Set $I=\ker(\nu)$ and $J=\ker(\mu)$, so that
$I=\varphi_{\alpha,u,\tau}(J)$. Suppose that $\alpha\sim_J u$ and let
$p\colon\pi_1(Q,I)\twoheadrightarrow \pi_1(Q,J)$ be the quotient
map (see Proposition~\ref{prop1.3}). Then, the following diagram
commutes:
\begin{equation}
\xymatrix{
\hom(\pi_1(Q,J),k^+) \ar@{->}[rd]^{\theta_{\mu}} 
\ar@{->}[dd]_{p^*}
& \\
& \hh^1(A)\\
\hom(\pi_1(Q,I),k^+) \ar@{->}[ru]_{\theta_{\nu}} & 
}\notag
\end{equation}
where $p^*\colon \hom(\pi_1(Q,J),k^+)\hookrightarrow
\hom(\pi_1(Q,I),k^+)$ is the embedding induced by $p$. In particular,
$\im(\theta_{\mu})\subseteq \im(\theta_{\nu})$.
\end{prop}
\noindent{\textbf{Proof:}} Recall that $p$ is the map
$[\gamma]_I\mapsto [\gamma]_J$.
Let $f\in \hom(\pi_1(Q,J),k^+)$. So
$p^*(f)$ is the composition
$\pi_1(Q,I)\xrightarrow{p}\pi_1(Q,J)\xrightarrow{f} k$. We know that
$\theta_{\mu}(f)$ and $\theta_{\nu}(p^*(f))$ are represented by the
derivations $d_1$ and $d_2$ respectively, such that for any
arrow $x\xrightarrow{a}y$:
\begin{align}
&d_1(\mu(a))=f([\gamma_y^{-1}a\gamma_x]_J)\
\mu(a)=p^*(f)([\gamma_y^{-1}a\gamma_x]_I)\ \mu(a)\notag\\
&d_2(\nu(a))=p^*(f)([\gamma_y^{-1}a\gamma_x]_I)\ \nu(a)\ .\notag
\end{align}
Let us prove that $d_1$ and $d_2$ coincide on $\nu(Q_1)$. Let
$x\xrightarrow{a}y$ be an arrow. If $a\neq \alpha$, then
$\mu(a)=\nu(a)$ and the above characterizations of $d_1$ and $d_2$
imply that $d_1(\nu(a))=d_1(\mu(a))=d_1(\nu(a))$. Now assume that
$a=\alpha$ so that:
$\nu(a)=\mu(a)-\tau\mu(u)$ and
$[\gamma_y^{-1}a\gamma_x]_J=[\gamma_y^{-1}u\gamma_x]_J$ (recall that $a=\alpha\sim_Ju$).
Thus:
\begin{equation}
\begin{array}{rclc}
d_1(\nu(a)) &=&d_1(\mu(\alpha))-\tau\ d_1(\mu(u))\\
&=&f([\gamma_y^{-1}\alpha\gamma_x]_J)\ \mu(\alpha)-\tau
f([\gamma_y^{-1}u\gamma_x]_J)\ \mu(u)\\
&=&f([\gamma_y^{-1}\alpha\gamma_x]_J)\ (\mu(\alpha)-\tau\ \mu(u))\\
&=&p^*(f)([\gamma_y^{-1}\alpha\gamma_x]_I)\ \nu(\alpha)\\
&=&d_2(\nu(\alpha))=d_2(\nu(a))&.
\end{array}\notag
\end{equation}
Hence, $d_1$ and $d_2$ are two derivations of $A$ and they coincide on
$\nu(Q_1)$. So $d_1=d_2$ and
$\theta_{\mu}(f)=\theta_{\nu}(p^*(f))$ for any $f\in
\hom(\pi_1(Q,J),k^+)$.\hfill$\square$\\

The following example shows that
Proposition~\ref{prop3.2} does not necessarily hold true if $\nu$ is a
presentation of $A$ and $\psi\colon kQ\to kQ$ is an automorphism such
that the identity map on the walks in $Q$ induces a surjective
group homomorphism $\pi_1(Q,\ker(\nu))\twoheadrightarrow\pi_1(Q,\ker(\nu\circ\psi))$.
\begin{ex}
\label{ex3.5}
Let $A=kQ/I$ where $\carac(k)=2$, where $Q$ is the quiver of Example~\ref{ex3.4} and
where $I=<da,fea+dcb>$. Let $\nu\colon kQ\twoheadrightarrow A$ be the
natural projection with kernel $I$, let
$\psi:=\varphi_{d,ef,1}\varphi_{a,cb,1}$ and let
$\mu:=\nu\circ\psi\colon kQ\twoheadrightarrow A$. Hence
$\ker(\mu)=<da+fecb, fea+dcb>$. Note that
$\pi_1(Q,\ker(\nu))\simeq\mathbb{Z}$ is generated by $[b^{-1}c^{-1}a]_I$ and
that $\pi_1(Q,\ker(\mu))\simeq\mathbb{Z}/2\mathbb{Z}$ is generated by
$[b^{-1}c^{-1}a]_J$. Note also that $\sim_{\ker(\nu)}$ is weaker that $\sim_{\ker(\mu)}$
so that the identity map on the set of walks in $Q$ induces a
surjective group homomorphism
$p\colon \pi_1(Q,\ker(\nu))\twoheadrightarrow\pi_1(Q,\ker(\mu))$. Let $T$ be the
maximal tree such that $T_1=\{b,c,e,f\}$. Let $f\colon
\pi_1(Q,\ker(\mu))\to k$ be the group homomorphism such that
$f([b^{-1}c^{-1}a]_J)=1$. On the one hand, $\theta_{\mu}(f)\in \hh^1(A)$
is represented by the derivation:
\begin{equation}
\begin{array}{rrclc}
d_1\colon & A & \rightarrow & A\\
& \mu(x) & \mapsto & \mu(x)\ \text{if $x\in\{a,d\}$}\\
& \mu(x) & \mapsto & 0\ \text{if $x\in \{b,c,e,f\}$}&.
\end{array}\notag
\end{equation}
On the other hand, $\theta_{\nu}(p^*(f))\in \hh^1(A)$ is represented by
the derivation:
\begin{equation}
\begin{array}{rrclc}
d_2\colon & A & \rightarrow & A\\
& \mu(a) & \mapsto & \mu(a)+\mu(cb)\\
& \mu(d) & \mapsto & \mu(d)+\mu(fe)\\
& \mu(x) & \mapsto & 0\ \text{if $x\in \{b,c,e,f\}$}&.
\end{array}\notag
\end{equation}
One checks easily that $d_2-d_1$ is not inner so that $\theta_{\mu}(f)\neq\theta_{\nu}(p^*(f))$.
Moreover, $\im(\theta_{\nu})$ and $\im(\theta_{\mu})$ are one
dimensional (because $\carac(k)=2$, $\pi_1(Q,\ker(\nu))\simeq\mathbb{Z}$
and $\pi_1(Q,\ker(\mu))\simeq\mathbb{Z}/2\mathbb{Z}$) and
$d_1,d_2$ are not inner. Hence $\im(\theta_{\mu})\not\subseteq \im(\theta_{\nu})$.

Actually, Proposition~\ref{prop3.2} does not work here because
the automorphism $\psi\colon (kQ,J)\to (kQ,I)$
maps arrows to linear combination of paths which are not homotopic
for $\sim_I$. For example, $\psi(a)=a+cb$ whereas
$a\not\sim_Icb$ (recall that
$\pi_1(Q,I)\simeq \mathbb{Z}$ is generated by $[b^{-1}c^{-1}a]_I$).
\end{ex}

Finally, we compare $\theta_{\nu}$ and $\theta_{\mu}$ when
$\mu=\nu\circ \psi$ with $\psi\colon kQ\xrightarrow{\sim}kQ$ an automorphism such that
$\ker(\nu)=\ker(\mu)$.
\begin{prop}
\label{prop3.3}
Let $\nu\colon kQ\twoheadrightarrow A$ be a presentation and let
$I=\ker(\nu)$. Let $\psi\colon kQ\xrightarrow{\sim} kQ$ be an
automorphism such that $\psi(e_i)=e_i$ for every $i$ and such that
$\psi(I)=I$. Let $\mu:=\nu\circ \psi\colon kQ\twoheadrightarrow A$ so
that $\ker(\mu)=I$. Let $\overline{\psi}\colon A\xrightarrow{\sim}A$ be the
$k$-algebra automorphism such that $\overline{\psi}\circ\mu=\mu\circ\psi$.
Then, the following diagram commutes:
\begin{equation}
\xymatrix{
& \hh^1(A) \ar@{->}[dd]^{\overline{\psi}_*}_{\sim}\\
\hom(\pi_1(Q,I),k^+) \ar@{->}[ru]^{\theta_{\nu}}
\ar@{->}[rd]_{\theta_{\mu}} & \\
& \hh^1(A)&.
}\notag
\end{equation}
In particular, $\im(\theta_{\mu})$ is equal to the image of
$\im(\theta_{\nu})$ under the Lie algebra automorphism
$\overline{\psi}_*\colon \hh^1(A)\xrightarrow{\sim} \hh^1(A)$ induced by
$\overline{\psi}\colon A\xrightarrow{\sim} A$.
\end{prop}
\noindent{\textbf{Proof:}} Since $\psi$ fixes the idempotents
$e_1,\ldots,e_n$, we know that $\overline{\psi}_*$ is well
defined. Let $f\in \hom(\pi_1(Q,I),k^+)$. So $\theta_{\nu}(f)$ and
$\theta_{\mu}(f)$ are represented  by the derivations $d_1$ and $d_2$
respectively, such that for any arrow $x\xrightarrow{\alpha}y$:
\begin{align}
&d_1(\nu(\alpha))=f([\gamma_y^{-1}\alpha\gamma_x]_I)\
  \nu(\alpha)\notag\\
&d_2(\mu(\alpha))=f([\gamma_y^{-1}\alpha\gamma_x]_I)\
  \mu(\alpha)\ .\notag
\end{align}
In order to prove that $\overline{\psi}_*(\theta_{\nu}(f))=\theta_{\mu}(f)$
it suffices to prove that $\overline{\psi}\circ
d_1=d_2\circ\overline{\psi}$. Let $x\xrightarrow{\alpha}y$ be an
arrow. Then:
\begin{equation}
\begin{array}{rcll}
d_2\circ\overline{\psi}(\nu(\alpha))&=&d_2(\mu(\alpha))&\text{because
  $\nu=\mu\circ\psi^{-1}$ and
  $\overline{\psi}\circ\mu=\mu\circ\psi$}\\
&=&f([\gamma_y^{-1}\alpha\gamma_x]_I)\ \mu(\alpha)&.
\end{array}\notag
\end{equation}
On the other hand:
\begin{equation}
\begin{array}{rcll}
\overline{\psi}\circ
d_1(\nu(\alpha))&=&
f([\gamma_y^{-1}\alpha\gamma_x]_I)\ \overline{\psi}(\nu(\alpha))&\\
&=&f([\gamma_y^{-1}\alpha\gamma_x]_I)\ \mu(\alpha) & \text{because
  $\mu=\nu\circ\psi$ and $\overline{\psi}\circ\mu=\mu\circ\psi$.}
\end{array}\notag
\end{equation}
Hence, $\overline{\psi}\circ d_1$ and $d_2\circ \overline{\psi}$ are derivations of
$A$ which coincide on $\nu(Q_1)$. So
$\overline{\psi}_*(\theta_{\nu}(f))=\theta_{\mu}(f)$ for any $f\in
\hom(\pi_1(Q,I,k^+))$.\hfill$\square$\\

\section{Proof of Theorem~\ref{thm1}}
In this section, we prove Theorem~\ref{thm1}. We begin with the
following useful lemma.
\begin{lem}
\label{lem4.1}
Assume that at least one of the two conditions ($H_1$) or ($H_2$) is satisfied.
Let $\nu\colon kQ\twoheadrightarrow A$ be a presentation whose kernel
$I_0$ is such that $\sim_{I_0}$ is the
unique source of $\Gamma$. Let $\mu\colon kQ\twoheadrightarrow A$
be another presentation. Then, there exist $\nu'\colon kQ\twoheadrightarrow
A$ a presentation with kernel $I_0$ and a $k$-algebra automorphism
$\psi\colon A\xrightarrow{\sim}A$ such that:
\begin{enumerate}
\item[-] $\psi(e_i)=e_i$ for any $i$,
\item[-] $\im(\theta_{\mu})\subseteq
\im(\theta_{\nu'})=\psi_*(\im(\theta_{\nu}))$.
\end{enumerate}
If moreover $\sim_I$ and $\sim_{I_0}$ coincide, then the above inclusion is an equality.
\end{lem}
\noindent{\textbf{Proof:}}
Let $\overline{\nu}\colon kQ/I_0\xrightarrow{\sim} A$ and $\overline{\mu}\colon
kQ/I\xrightarrow{\sim} A$ be the isomorphisms induced by $\nu$
and $\mu$ respectively. Hence, $\overline{\mu}^{-1}\circ\overline{\nu}\colon
kQ/I_0\xrightarrow{\sim} kQ/I$ is an isomorphism which maps $e_i$ to
$e_i$ for every $i$. Hence, there exists an automorphism
$\varphi\colon kQ\xrightarrow{\sim} kQ$ which maps $e_i$ to $e_i$ for
every $i$ and such that the following diagram commutes (see 
\cite[Prop. 2.3.18]{lemeur_thesis}, for instance):
\begin{equation}
\xymatrix{
kQ \ar@{->>}[d] \ar@{->}[r]^{\varphi}& kQ \ar@{->>}[d]\\
kQ/I_0\ar@{->}[r]^{\bar{\mu}^{-1}\bar{\nu}}  & kQ/I&.}\notag
\end{equation}
where the vertical arrows are the natural projections. So,
the following diagram is commutative:
 \begin{equation}
 \xymatrix{
 kQ\ar@{->}[rr]^{\varphi}
\ar@{->}[rd]_{\nu}
&&
kQ\ar@{->}[ld]^{\mu}
\\
 &A&}\notag
 \end{equation}
Let us apply Proposition~\ref{prop1.4} to $I$. 
We keep the notations $\alpha_i,\,u_i,\,\tau_i,\,I_i$ of that
proposition. Let
$\psi:=\varphi^{-1}D\varphi_{\alpha_l,u_l,\tau_l}\ldots\varphi_{\alpha_1,u_1,\tau_1}$.
Thus, $\psi(I_0)=I_0$ and 
$\nu':=\nu\circ\psi\colon kQ\twoheadrightarrow A$ is a presentation
with kernel $I_0$.
Thanks to Proposition~\ref{prop3.3} from which we keep the notations,
we know that:
\begin{equation}
\im(\theta_{\nu'})=\overline{\psi}_{\star}(\im(\theta_{\nu}))\ .\tag{$1$}
\end{equation}
Let us show that $\im(\theta_{\mu})\subseteq
\im(\theta_{\nu'})$. By construction, we have
$\mu=\nu\varphi^{-1} =
\nu'\varphi_{\alpha_1,u_1,\tau_1}^{-1}\ldots\varphi_{\alpha_l,u_l,\tau_l}^{-1}D^{-1}$.
For simplicity, we  use the following notations: 
$\mu_0:=\nu'$ and
$\mu_i:=\nu'\varphi_{\alpha_1,u_1,\tau_1}^{-1}\ldots\varphi_{\alpha_i,u_i,\tau_i}^{-1}$
for $i\in\{1,\ldots,l\}$.
Note that $I_i=\ker(\mu_i)$. Since $\alpha_i\sim_{I_i}u_i$ and
$\mu_i=\mu_{i-1}\circ\varphi_{\alpha_i,u_i,-\tau_i}$,
Proposition~\ref{prop3.2} implies that:
\begin{equation}
\im(\theta_{\mu_l})\subseteq \im(\theta_{\mu_{l-1}})\subseteq \ldots
\subseteq \im(\theta_{\mu_i})\subseteq
\im(\theta_{\mu_{i-1}})\subseteq\ldots\subseteq
\im(\theta_{\mu_0})=\im(\theta_{\nu'})\ .\tag{$2$}
\end{equation}
Moreover,
 $\mu=\mu_mD^{-1}$ where
 $D^{-1}$ is a dilatation. Hence, ($1$), ($2$) and Proposition~\ref{prop3.1} imply that:
\begin{equation}
\im(\theta_{\mu})=\im(\theta_{\mu_l})\subseteq
\im(\theta_{\nu'})=\overline{\psi}_{\star}(\im(\theta_{\nu}))\ .\tag{$3$} 
\end{equation}
Now assume that $\sim_I$ is the unique source of $\Gamma$. Then
Proposition~\ref{prop1.4} imply that the homotopy relations
$\sim_{I_0},\sim_{I_1},\ldots,\sim_{I_l},\sim_I$ coincide.
Therefore, for any $i\in\{1,\ldots,l\}$, we have
$\mu_{i-1}=\mu_i\circ\varphi_{\alpha_i,u_i,\tau_i}$,
and 
$\alpha_i\sim_{I_{i-1}}u_i$.
So Proposition~\ref{prop3.2}, implies that
$\im(\theta_{\mu_{i-1}})\subseteq \im(\theta_{\mu_i})$. This proves that
all the inclusions in $(2)$ are equalities, and so is the inclusion
in $(3)$.\hfill$\square$\\

Now we can prove Theorem~\ref{thm1}.\\
\noindent{\textbf{Proof of Theorem~\ref{thm1}:}}
$(i)$ Let $\mathcal{G}$ be a 
maximal diagonalizable subalgebra of $\hh^1(A)$. Thanks to
Proposition~\ref{prop2.6}, there exists a presentation
$\mu\colon kQ\twoheadrightarrow A$ such that $\mathcal{G}\subseteq
\im(\theta_{\mu})$. On the other hand, Lemma~\ref{lem4.1}, implies that
there exists a presentation $\nu\colon kQ\twoheadrightarrow A$
such that $\sim_{\ker(\nu)}$ is the unique source of $\Gamma$
and such that $\im(\theta_{\mu})\subseteq \im(\theta_{\nu})$.
Hence, $\mathcal{G}\subseteq \im(\theta_{\nu})$ where
$\im(\theta_{\nu})$ is a diagonalizable subalgebra of $\hh^1(A)$, thanks
to Proposition~\ref{prop2.4}. The maximality of
$\mathcal{G}$ forces $\mathcal{G}=\im(\theta_{\nu})$.

Conversely, let $\mu\colon kQ\twoheadrightarrow A$ be a presentation
such that $\sim_{\ker(\mu)}$ is the unique source of
$\Gamma$. Hence, $\im(\theta_{\mu})$ is diagonalizable (thanks to
Proposition~\ref{prop2.4}) so there exists a maximal
diagonalizable subalgebra $\mathcal{G}$ of $\hh^1(A)$ containing
$\im(\theta_{\mu})$. Thanks to the above description, we know that
$\mathcal{G}=\im(\theta_{\nu})$ where $\nu\colon
kQ\twoheadrightarrow A$ is a presentation such that $\sim_{\ker(\nu)}$ is the unique
source of $\Gamma$. 
Moreover, Lemma~\ref{lem4.1} gives a $k$-algebra automorphism $\psi\colon A\xrightarrow{\sim} A$
such that $\im(\theta_{\mu})=\psi_{\star}(\im(\theta_{\nu}))$. Since
$\psi_*$ is a Lie algebra automorphism of $\hh^1(A)$, the maximality of
$\mathcal{G}=\im(\theta_{\nu})$ implies that $\im(\theta_{\mu})$ is maximal.

$(ii)$ is a consequence of $(i)$ and of  Lemma~\ref{lem4.1}.\hfill$\square$\\

\bibliographystyle{plain}
\bibliography{biblio}

\noindent Patrick Le Meur\\
\textit{e-mail:} Patrick.LeMeur@cmla.ens-cachan.fr\\
\textit{address:} CMLA, ENS Cachan, CNRS, UniverSud, 61 Avenue du President Wilson, F-94230 Cachan
\end{document}